\documentclass[11pt,leqno,english]{article}
 \usepackage{latexsym}
 \usepackage{amssymb}
 \usepackage{amsmath}
 \usepackage{amsthm}
\usepackage{natbib}

 \usepackage[dvips]{graphics}
\usepackage{graphicx}

 \newcommand{ \un }{\mathbb{I} }
 \newcommand{ \p }{\mathbb{P} }
 \newcommand{ \pa }{\mathbb{P}^{\alpha} }

 \newcommand{ \E }{\mathbb{E}}

 \newcommand{ \R }{ \mathbb{R} }
 \newcommand{ \Z }{ \mathbb{Z} }
 \newcommand{\N}{ \mathbb{N} }
 
 \newcommand{ \V }{\textrm{Var} }

 \newcommand{ \f }{ \mathcal{F} }

 \newcommand{ \tm }{ m }
 \newcommand{ \tM }{ M }

\newcommand{ \peq }{\dot{=} }

\newtheorem{The}{{\bf Theorem}}[section]
\theoremstyle{definition}
\newtheorem{Def}[The]{{\bf Definition}}
\theoremstyle{plain}

 \newtheorem{Pro}[The]{\bf Proposition}
 \theoremstyle{definition}
\newtheorem{Rem}[The]{{\bf Remark}}

 \newenvironment{Pre}{\noindent \textbf{Proof.} \\ }{$\
 \blacksquare$}

\makeatletter

\@addtoreset{equation}{section} \makeatother

\title{A limit result for a system of particles in random environment  
\\ \vspace{1cm}
 \large{Pierre Andreoletti} \footnote{ Laboratoire MAPMO - C.N.R.S. UMR 6628 - F\'ed\'eration Denis-Poisson, Universit\'e d'Orl\'eans, 
(Orl\'eans France). Supported  by GREFI-MEFI and Departimento di Mathematica, Universita di Roma II ''Tor Vergata'', Italy. \newline \vspace{0.1cm}. \newline  MSC 2000 60G50; 82B41; 82B21. \newline \vspace{0.5cm} \textit{Key words and phrases :  random environment, random walk,
Sinai's regime, system of particles, excursion.} }   }

\date{August 28, 2007}

\begin{document}

\bibliographystyle{unsrtnat}
\maketitle


\noindent \\ \textbf{Abstract:} We consider an infinite system of particles in one dimension, each particle performs independant Sinai's random walk in
random environment. Considering an instant $t$, large enough, we prove a result in probability showing that the particles are trapped in the neighborhood of well defined points of the lattice depending on the random environment the time $t$ and the starting point of the particles.

\section{Introduction, definitions and result}
Systems of particles have been study a lot in many different directions see for example \cite{Ligget}, \cite{DeMPre} and \cite{KipLan} and the references therein.
In this paper we are interested in a system of particles performing independant Sinai's random walk in random environment on the latice $\Z$. 
At time 0, in each point of the lattice $i\in \Z$ a random number of particles $\eta(i,0)$ is distributed, we assume that the sequence $(\eta(i,0),i\in \Z)$ is i.i.d. and that $\eta(0,0)$  is distributed according to a Poisson distribution of parameter $\lambda$. 
Each particle of this system effectuates a Sinai's walk that means a one
dimensional random walk in random environment with three
conditions on the random environment: two necessaries hypothesis
to get a recurrent process (see \cite{Solomon}) which is not a
simple random walk and an hypothesis of regularity which allows us
to have a good control on the fluctuations of the random
environment. 
The asymptotic behavior of  such a walk has been understood by
\cite{Sinai} : it is sub-diffusive, that is, for an instant $t$ it behaves like $(\log t)^2$ and it is localized in the neighborhood of a well
defined point of the lattice.
We are interested in the general behavior of the system described above in the following sense, let us denote $\eta(x,t)$ the number of particles at time $t$ on $x$, and let $f$ be a differentiable function with compact support, we study the asymptotic behavior when $t$ becomes large of
\begin{eqnarray}
 \frac{1}{(\log t)^2}\left(\sum_{x\in \Z}\eta(x,t)f(x/(\log t)^2)\right).
\end{eqnarray}
The result we get can been seen as the generalization of Sinai's localization result for a system of particles. It shows that the particles are trapped in the bottom of the deepest valley that the walk can reach within an amount of time $t$. The coordinates of the bottom of the valleys depend on the starting points of the particles, the random environment and $t$. 

To get the result, first, we make  a construction on the random environment closed to the one that \cite{NevPit}  use to study the excursion of the Brownian motion, thanks to this we get a partition in term of valleys of the part of the lattice  related to the support of $f$. Then we prove that there is neither importation of particles from outside of the support of the function $f$ nor  exchanged of particles between the (deepest) valleys, finally we  we use the accurate result of \cite{Pierre1} about Sinai's localization result extended to a finite number of particles.

In the following section, we define Sinai's walk, then the system of particles and finally, we state the main result.
In section 2 we make the construction on the random environment, in section 3 and 4 we prove the main theorem. 

\subsection{Definition of Sinai's walk}

Let $\alpha =(\alpha_i,i\in \Z)$ be a sequence of independent and identically distributed (i.i.d.) random
variables taking values in $(0,1)$ defined on the probability
space $(\Omega_1,\f_1,Q)$, this sequence will be called random
environment. A random walk in random environment (denoted
R.W.R.E.) $(X_n,n
\in\N)$ is a sequence of random variable taking value in $\Z$, defined on $( \Omega,\f,\p)$ such that \\
$ \bullet $ for every fixed environment $\alpha$, $(X_n,n\in \N)$
 is a Markov chain with the following transition probabilities, for
 all $n\geq 1$ and $i\in \Z$
 \begin{eqnarray}
& & \p^{\alpha,1}\left[X_n=i+1|X_{n-1}=i\right]=\alpha_i, \label{mt} \\
& & \p^{\alpha,1}\left[X_n=i-1|X_{n-1}=i\right]=1-\alpha_i. \nonumber
\end{eqnarray}
We denote $(\Omega_2,\f_2,\p^{\alpha,1})$ the probability space
associated to this Markov chain. \\
 $\bullet$ $\Omega = \Omega_1 \times \Omega_2$, $\forall A_1 \in \f_1$ and $\forall A_2 \in \f_2$,
$\p\left[A_1\times
A_2\right]=\int_{A_1}Q(dw_1)\int_{A_2}\p^{\alpha(w_1),1}(dw_2)$.

\noindent \\ The probability measure $\p^{\alpha,1}\left[\left.
.\right|X_0=a \right]$  will be  denoted $\p^{\alpha,1}_a\left[.\right]$,
 the expectation associated to $\p^{\alpha,1}_a$: $\E^{\alpha,1}_a$, and the expectation associated to $Q$:
 $\E_Q$. 

\noindent \\ Now we introduce the hypothesis we will use in all
this work. The two following hypothesis are the necessaries
hypothesis
\begin{eqnarray}
 \E_Q\left[ \log
\frac{1-\alpha_0}{\alpha_0}\right]=0 , \label{hyp1bb} \label{hyp1}
\end{eqnarray}
\begin{eqnarray}
\V_Q\left[ \log \frac{1-\alpha_0}{\alpha_0}\right]\equiv \sigma^2
>0 . \label{hyp0}
\end{eqnarray}
 \cite{Solomon} shows that under \ref{hyp1} the
process $(X_n,n\in \N)$ is $\p$ almost surely recurrent and
\ref{hyp0} implies that the model is not reduced to the simple
random walk. In addition to \ref{hyp1} and \ref{hyp0} we will
consider the following hypothesis of regularity, there exists $0<
\rho_0 < 1/2$ such that
\begin{eqnarray}
& & \sup \left\{x,\ Q\left[\alpha_0 \geq x \right]=1\right\} \geq \rho_{0} \textrm{ and } \sup
\left\{x,\ Q\left[\alpha_0 \leq 1-x \right]=1\right\} \geq \rho_0.
\label{hyp4}
\end{eqnarray}

\noindent We call \textit{Sinai's random walk} the random walk in
random environment previously defined with the three hypothesis
\ref{hyp1}, \ref{hyp0} and \ref{hyp4}.
Instead of considering the discrete time process, we will work with the continuous time one $(X_t,t\in \R_+)$ define as follows, $\forall t \in \R_+,\ X_t=X_{[t]}$ where $[t]$ is the integer par of $t$.

\subsection{The system of particles}

We recall that $\eta(i)\equiv\eta(i,0)$ is the number of particles  in $i\in \Z$ at time $0$, we assume that   $(\eta(i),i\in \Z)$ is i.i.d, and that $\eta(0)$ is distributed according to a Poisson distribution of parameter $\lambda$, we denote $\p_{1}$ the corresponding probability measure. We denote $X_{t}^{x,i}$ the coordinate of a particle at time $t$  corresponding to the $i^{\textrm{th}}$ particle,  which was on site $x$ at time $0$ $(X^{x,i}_{0}=x)$.
\noindent $\eta(x,t)$ can be written in the following useful way: 
\begin{eqnarray}
 \eta(x,t)  & = &\sum_{y \in \Z} \sum_{i=1}^{\eta(y)}  \un_{\{X_{t}^{x,i}=x\}}
 \end{eqnarray}
The distribution of the whole system is denoted $P$, it is the product of $Q$ and $\pa$, where $\pa$ is the product $\pa \equiv  \prod_{x \in \Z}\p_x^{\alpha}=P_1\times \prod_{x \in \Z}\prod_{i=1}^{\eta(x)} \p_x^{\alpha,i}$. The measure $P$ generalizes the measure $\p$ defined Section 1.1 for a random number of particles. Here we have three levels of randomness, the environment, the number of particles per sites and the random walks themselves, hence both $\pa$ and $\pa_x$ are random measures. 
Notice that, when the environment is fixed $(\eta(.,t))$ inherits the markov property from the random walks, indeed they all have the same transition probability and they do not interact with each other.
 \subsection{Main results}

First let us state the result in probability, the following Theorem is based on a construction on the random environment (see figure \ref{fig8}). Let $C^1_{\kappa}$ be the set of all differentiable functions with a compact support. 
 

\begin{The} \label{th1}
Assume \ref{hyp1bb}, \ref{hyp0} and \ref{hyp4} hold, for all $f\in C^1_{\kappa}$, and all $t$ large enough
\begin{eqnarray}
 & & P \left[  \left| \frac{1}{(\log t)^2} \sum_{x\in \Z}\eta(x,t)f(x/(\log n)^2)-  \lambda \sum_{i=1}^{n(f)}\frac{|M_{i+1}-M_{i}|}{(\log t)^2}f(m_{i}/(\log t)^2) \right| = o(1) \right] = 1- o(1),  \nonumber
\end{eqnarray}
where, $(m_{i},i\geq 1)$, $(M_{i},i \geq 1)$ and $n(f)$  are well defined variables depending only on the random environment $\alpha$, on the time $t$ and $f$. $o(1)$ is a positive decreasing function such that $\lim_{t \rightarrow + \infty}o(1)=0$.
\end{The}

\noindent This Theorem says that at time $t$ the particles are localized in the neighborhood of  well defined points of the lattice. These points have the coordinates of the bottoms of the valleys where the particles started their walk. 
In other words it says that the number of particles at time $t$, at the points $m_i$ is equal to $\lambda |M_{i+1}-M_{i}|$ which is the total number of particles at the instant $t=0$ present in the valley $(\{M_{i},m_{i},M_{i+1}\})$. $n(f)$ is the number of valleys deep enough within the support of $f$,  notice that $n(f)$ can be equal to zero if the support of $f$ is too small. The following Proposition give the limit distributions of the random variables present in the Theorem. 

\begin{Pro} \label{ch1}
Assume \ref{hyp1bb}, \ref{hyp0} and \ref{hyp4} hold, for all $i \in \N^*$, in distribution when $t$ goes to infinity
\begin{eqnarray}
\frac{\sigma^2M_{i}}{(\log t)^2} \rightarrow L_i, \\
\frac{\sigma^2m_{i}}{(\log t)^2} \rightarrow L_i,
\end{eqnarray}
The $(L_{i+1}-L_i,i \in \Z )$ are independent and equidistributed random variables with Laplace transform given by $E(\exp(-\lambda (L_{2}-L_1))) = 1/(\cosh(\sqrt( {2} \lambda)))^2$, for all $\lambda>0$. 
Moreover when the support of $f$ is large enough the distribution of $n(f)$ is given by a normal law with mean $\sigma^2 supp(f)/2$ and variance $3* \sigma^4 supp(f)/4$.
\end{Pro}

\noindent


\section{ Construction for the random environment \label{sec11}}
In this section we begin with some basic notions on the random environment of Sinai's random walk, then 
we construct the key random variables $(M_{i},m_{i}, i \leq n(f))$ appearing in our result. The reader can follow the different steps of the construction with the Figures 2-5.

\subsection{ Definition and basic notions of valleys \label{BasicS} }

For completeness we begin with some basic notions originally introduced by \cite{Sinai}.
\noindent \\ \textbf{The random potential and the valleys}
\noindent \\ Let
\begin{eqnarray}
\epsilon_i \equiv \log \frac{1-\alpha_i}{\alpha_i},\ i\in \Z,
\end{eqnarray}
define :
\begin{Def} \label{defpot2} The random potential $(S_k,\  k \in
\R)$ associated to the random environment $\alpha$ is defined in the following way: 
\begin{eqnarray} 
 && S_{k}=\left\{ \begin{array}{ll} \sum_{1\leq i \leq k} \epsilon_i, & \textrm{ if }\ k =1,2, \cdots, \\
 - \sum_{k+1 \leq i \leq 0} \epsilon_i , &   \textrm{ if }\  k=-1,-2,\cdots,   \end{array} \right.  \nonumber \\
&& S_{0}=0. \nonumber
 \end{eqnarray}
 for the other $k  \in \R \setminus \Z$, $S_k$ is defined by linear interpolation.
\end{Def}


\begin{Def} \label{c2s2d1}
 We will say that the triplet $\{M',m,M''\}$ is a \textit{valley} if
 \begin{eqnarray}
& &  S_{M'}=\max_{M' \leq t \leq m} S_t ,  \\
& &  S_{M''}=\max_{m \leq t \leq
 \tilde{M''}}S_t ,\\
& & S_{m}=\min_{M' \leq t \leq M''}S_t \ \label{2eq58}.
 \end{eqnarray}
If $m$ is not unique we choose the one with the smallest absolute
value.
 \end{Def}

\begin{Def} \label{deprofvalb}
 We will call \textit{depth of the valley} $\{\tM',\tm,\tM''\}$ and we
 will denote it $d([M',M''])$ the quantity
\begin{eqnarray}
 \min(S_{M'}-S_{m},S_{M''}-S_{m})
 .
 \end{eqnarray}
 \end{Def}

\noindent  Now we define the operation of \textit{refinement}
 \begin{Def} \label{Defref}
Let  $\{\tM',\tm,\tM''\}$ be a valley and let
  $\tM_1$ and $\tm_1$ be such that $\tm \leq \tM_1< \tm_1 \leq \tM''$
  and
 \begin{eqnarray}
 S_{\tM_1}-S_{\tm_1}=\max_{\tm \leq t' \leq t'' \leq
 \tM''}(S_{t'}-S_{t''}) .
 \end{eqnarray}
 We say that the couple $(\tm_1,\tM_1)$ is obtained by a \textit{right refinement} of $\{\tM',\tm,\tM''\}$. If the couple $(\tm_1,\tM_1)$ is not
 unique, we will take  the one such that $\tm_1$ and $\tM_1$ have the smallest  absolute value. In a similar way we
  define the \textit{left refinement} operation.
 \end{Def}

\noindent \\  We denote $\log_2\equiv\log \log $, in all this work we will suppose that $t$ is
large enough such that $\log_2 t$ is positive.

\begin{Def} \label{thdefval1b} Let $\gamma>0$, define $\Gamma_t \equiv \log t+ \gamma \log_2 t $, we say that a valley $\{\tM',\tm,\tM''\}$
is of depth larger than $\Gamma_t$ if and only if $d\left([\tM',\tM'']\right)\geq \Gamma_t $. 
\end{Def}

\subsection{ Construction of a random cover for $(supp(f))*(\log t)^2$  }

We start with the constructions of a sequence of valleys with identical properties, the different steps of the construction are shown in Figures \ref{fig6} to \ref{fig7}. Notice that the variables in greek letters we use below are not the fundamental one, they just help us to construct $m^+_{.}, m^-_{.}, M^+_{.}$, $M^-_{.}$ and $M_{0}$.

 \noindent Let  $u,v \in \R$, $u<v$, define
\begin{eqnarray}
& & S_{u,v}^-=\min\{S_s,\ u \leq s \leq v \}, \\
& & S_{u,v}^+=\max\{S_s,\ u \leq s \leq v \}, \\
& & \tau^+_0  =  \inf  \{s >0,\ S_s-S_{0,s}^- \geq \Gamma_{t} \} ,  \\
& & m_{0}^+ =  \sup  \{ s < \tau^+_0 ,\ S_{s}=S_{0,\tau^+_{0}}^-  \}, \\
& & \sigma_{0}^+ =  \inf  \{ s > m^+_0 ,\ S_{m^{+}_{0},s}^+-S_s \geq \Gamma_{t}  \}, \\
& & \tau^-_0  =  \sup  \{s <0,\ S_s-S_{s,0}^- \geq \Gamma_{t} \} ,  \\
& & m_{0}^- =  \sup  \{ s > \tau^-_0 ,\ S_{s}=S_{0,\tau^+_{0}}^-  \}, \\
& & \sigma_{0}^- =  \sup  \{ s < m^-_0 ,\ S_{m^-_{0},s}^+-S_s \geq \Gamma_{t}  \}, \\
& & M_{0}= \inf  \{s >  m^-_{0},\ S_{s}=S_{m^-_{0},m^+_{0}}^+ \}. 
\end{eqnarray}

\begin{figure}[h]
\begin{center}
\input{fig7.pstex_t} 
\caption{Construction of the valleys 1/2} \label{fig6}
\end{center}
\end{figure}
\noindent
Now define recursively the following variables, let $i\geq 1$
\begin{eqnarray}
& & \tau^+_{i}  =  \inf  \{s >  \sigma^+_{i-1},\ S_{s}-S_{ \sigma^+_{i-1},s}^+ \geq \Gamma_{t} \} ,  \\
& & m_{i}^+ =  \sup  \{ s > \sigma^+_{i-1} ,\ S_{s}=S_{\sigma^+_{i-1},\tau^+_{i}}^+  \}, \\
& & M_{i}^+=  \inf  \{s > m^+_{i-1},\ S_{s}=S_{ m_{i-1}^{+},m_{i}^+}^+ \} , \\
& & \sigma_{i}^+ =  \inf  \{ s > m^+_i ,\ S_{m^{+}_{i},s}^+-S_s \geq \Gamma_{t}  \}, 
\end{eqnarray}
and in the same way on the left hand side of the origin.

\begin{figure}[h]
\begin{center}
\input{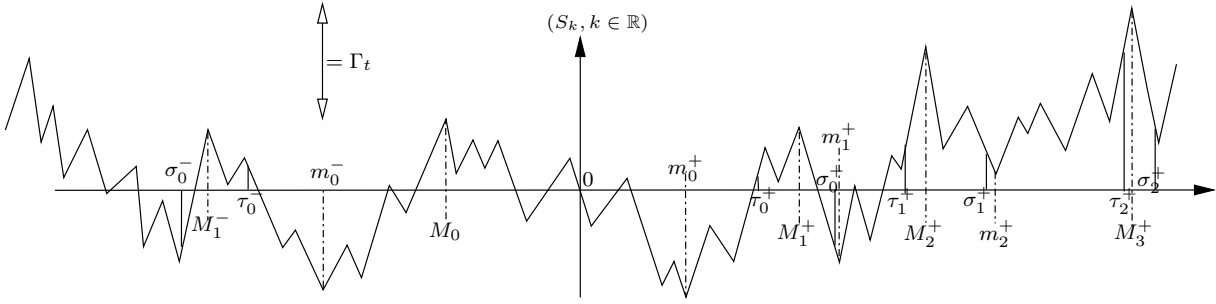} 
\caption{Construction of the valleys 2/2} \label{fig7}
\end{center}
\end{figure}

We would like to make remarks about this construction: 
\begin{Rem}
Notice that for all $i\geq 1$ the $\{M_{i}^+,m_{i}^+,M_{i+1}^+\}$ are valleys of depth larger than $\Gamma_{t}$ such that we can not construct a sub valley (with the operation of refinement, see Definition \ref{Defref} ) of depth larger than $\Gamma_{t}$, the same remark is true if we replace the $+$ by the $-$. 
However we do not know how deep the valleys $\{M_{1}^-,m_{0}^-,M_{0}\}$ and  $\{M_{0},m_{0}^+,M_{1}^+\}$ are, what we know from our construction is that we must have one of these two cases:
\begin{eqnarray}
& & \min\left (d\left([M_{1}^-,M_{0}]\right), d\left([M_{0},M_{1}^{+}]\right) \right)\geq \Gamma_{t} \textrm{ or} \label{case1}\\
& &  \min\left(d\left([M_{1}^-,M_{0}]\right), d\left([M_{0},M_{1}^{+}]\right) \right) < \Gamma_{t} \textrm{ and } d\left([M_{1}^-,M_{1}^{+}]\right) \geq \Gamma_{t} \label{case2}\
\end{eqnarray}
In our drawings it is the case \ref{case1} that occurs.
\end{Rem}

\noindent Now let us describe the link between the construction we have exposed above and the support of $f$. For simplicity, we will assume that $supp(f)=[-K,K]$ with $K>0$. Let $n^+$ be the largest integer such that $m_{n^+}<K (\log t)^2$ and symmetrically $n^-$ the largest integer such that $m_{n^-}>-K(\log t)^2$. Note that both  ${n^+}$ and ${n^-}$ can be equal to zero. It is also important to notice that, we can have 2 different cases in both on the right and left hand-side of the origin. We will only discuss about the right hand-side: first case is when $m_{n^+}<K(\log t)^2<M_{n^++1}$ in this case we will have to consider the particles that belong to the interval $(K(\log t)^2,M_{n^++1})$. Indeed most of them will tend to reach the bottom the valley $m_{n^+}$ that belongs to $[-K(\log t)^2, K(\log t)^2]$. To the contrary in the case where $M_{n(f)+1}>K(\log t)^2>m_{n^++1}$, we will not have to care about the particles in the interval $(M_{n(f)+1},K(\log t)^2)$ because they will all tend to reach $m_{n^++1}$ that do not belongs to  $[-K(\log t)^2, K(\log t)^2]$.

\begin{figure}[h]
\begin{center}
\input{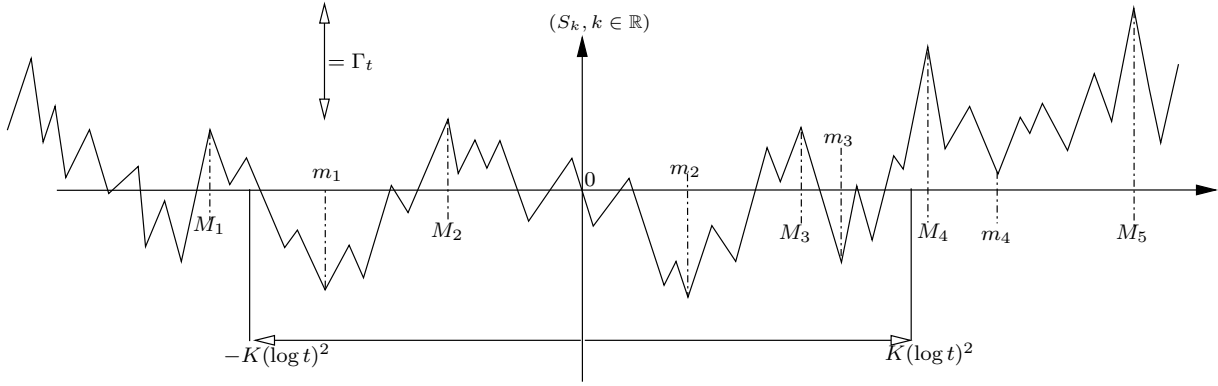} 
\caption{Support of $f$ and the valleys} \label{fig8}
\end{center}
\end{figure}

To simplify notations we re-numerate all the set of minima and maxima in the following way : we suppress the ''+'' and the ''-'', so denoting $n(f) \equiv n(K)=(n^++n^--\un_{A^c_t})\vee 0$ where $A_t=\{ \min\left (d\left([M_{1}^-,M_{0}]\right), d\left([M_{0},M_{1}^{+}]\right) \right)\geq \Gamma_{t}\}$, we denote: $M_{n^-}^-=M_1$, $m_{n^-}^-=m_1$, $M_{n^- -1}^-=M_{2}$, ..., $M_{n^+}^+=M_{n(f)+1}$. Note that the $-\un_{A^c_t}$ is needed because we have two cases (\ref{case1} and \ref{case2}). In our drawing (see Figure \ref{fig8}) $n^+=2$, $n^-=1$ and we are in the case \ref{case1} so $n(f)=3$.

\subsection{ Indeterminate points and Set of good environments}

Let $x \in [-K(\log t)^2,K(\log t)^2]$, due to the preceding construction we know that there exists $i\in \N$ such that $x \in [M_{i},M_{i+1}]$, let $\bar{i}\equiv \bar{i}(t,f,\alpha) =\{y \in \Z,\ y \in [M_{i},M_{i+1}] \}$, $i \leq n(f)$. 
Now let us define the set of \textit{indeterminate points}. Let $x \in \bar{l}$, $1 \leq l \leq n(f)+1$, we will say that $x$ is indeterminate if $x\in U_{l}\equiv [M_{l}-(\log_{2}t)^2,M_{l}+(\log_{2}t)^2] $. We call these sets "indeterminate", because at each top of the valley, actually a small interval on the top, we can not determine which valley the particles starting from there will choose. Notice, that the number of indeterminate points, belonging to $V_f$, are negligible comparing to a typical fluctuation of Sinai's random walk.

\begin{figure}[h]
\begin{center}
\input{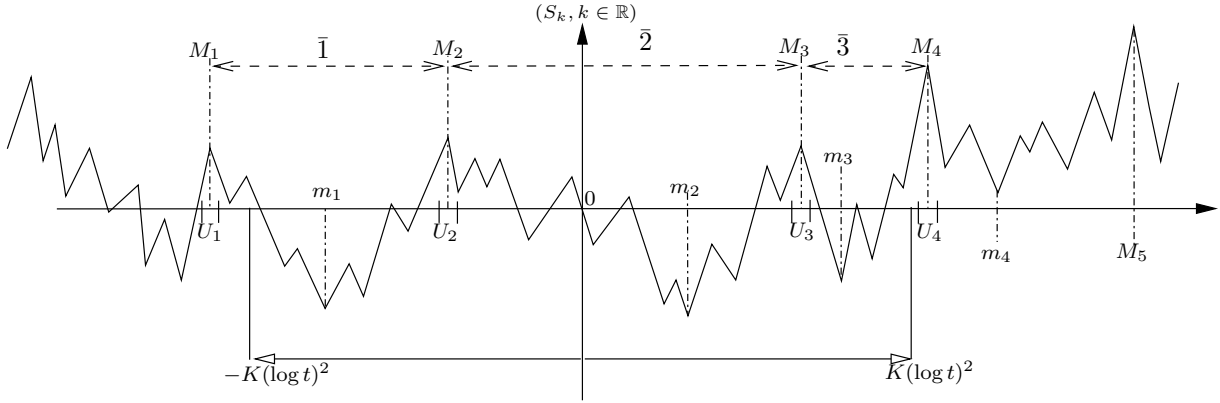} 
\caption{Random cover and sets of indeterminate points  with $n(f)=3$.} \label{fig9}
\end{center}
\end{figure}

\noindent
Now we give the main properties of the random environment, recall that $f\in C^1_{\kappa}$, $c_1>0$ and $c_2>0$:
\begin{eqnarray}
& & n(f) \textrm{ is finite}, \label{pr1}\\
& & |V_f|\leq c_1 (\log n)^2 \label{pr2}, \\ 
& & \forall 1 \leq i \leq n(f),\ |\bar{i}| \equiv |M_{i}-M_{i+1}| \thickapprox c_2 (\log n)^2 \label{pr3}, \\
&  & \forall 1 \leq i \leq n(f) \textrm{ the valley} \{M_{i},m_i,M_{i+1}\} \textrm{ have Sinai's basic properties} \label{pr4}.
\end{eqnarray}

\begin{Rem} (about \ref{pr4}) As we said before to get our result we need to prove the result of localization of Sinai for independent particles. For that we need some properties on the random environment we do not want to repeat here (see for example \cite{Pierre1} Definition 3.4, page 894), we have called these properties "Sinai's basic properties".
\end{Rem}
\noindent 
Let us define the set of good properties $G_{t}$, $G_{t}=\{\alpha \in \Omega_{1}, \alpha\ \textrm{ satisfies } \ref{pr1}-\ref{pr4} \}$. The following Proposition, give a probability result for $G_{t}$:

\begin{Pro}  \label{proenv}
There exists $c_{1}>0$, $c_2>0$ such that for all $t$ large enough, $Q[G_{t}] = 1- o(1)$.
\end{Pro}

\begin{Pre}
We shortly discuss the proof of this Proposition at the end of the paper, Paragraph 5. 
\end{Pre}

\section{Results of localization for the particles}
In this section we first prove two results (Proposition \ref{3.1}), about the no-migration of the particles from the valley where they start to another (different) one, then we state the key result which generalizes Sinai's result for a system of independent particles, let $f\in C ^1_{\kappa} $ we denote $V_{f}=\{\bar{1},\cdots,\bar{n}(f)\}$ the random cover of $[-K(\log t)^2, K(\log t)^2]$ (see the construction in the previous paragraph), we have
 \begin{Pro} \label{3.1} For all $f\in C^1_{\kappa}$, $\gamma>12$, and all $t$ large enough, $Q[G_{t}] = 1-o(1)$ and 
\begin{eqnarray}
\inf_{\alpha \in G_{t}}\left\{\pa\left[\sum_{x \notin V_{f}}  \sum_{k=1}^{\eta(x)}\un_{X_{t}^{x,k} \in V_{f}}=o(\log^2t) \right]\right\} = 1-o(1), \label{3.1eq1}
 \end{eqnarray}
conversly, we have:
\begin{eqnarray}
\inf_{\alpha \in G_{t}}\left\{\pa\left[\bigcap_{i=1}^{n(f)}\bigcap_{x\in 
\bar{i} \smallsetminus U_{i}}\bigcap_{k=1}^{\eta(x)} \left\{ X_{t}^{x,k} \in \bar{i} \right\}\right]\right\} = 1-o(1). \label{3.1eq2}
 \end{eqnarray}
 \end{Pro}
 
\noindent \ref{3.1eq1} shows that with an overwhelming probability a negligible number of particles, coming from the complementary of $V(f)$ in $\Z$, will reach  $V_{f}$ before the instant t.  \ref{3.1eq2} shows that, within the interval of time $t$, there is no exchange of particles between two valleys except, may be, for the particles that belong to the sets of indeterminate points.

\begin{Pre}
First let us give the proof of \ref{3.1eq2}, it is based on the two following facts,
\begin{eqnarray}
& &\pa\left[\bigcap_{i=1}^{n(f)}\bigcap_{x\in 
\bar{i} \smallsetminus U_{i}}\bigcap_{k=1}^{\eta(x)} \left\{ X_{t}^{x,k} \in \bar{i} \right\}\right] \geq \exp\left(\lambda \left[ a-1\right] \sum_{k=1}^{n(f)}|\bar{k}\smallsetminus U_k|\right),  \forall \alpha, \label{fact1} \\
&& a \geq 1- \frac{cte}{(\log t)^{\gamma-2}},\ \forall \alpha \in G_t. \label{fact2}
\end{eqnarray}
where $a=\inf_{1 \leq i \leq n(f)} \inf_{x \in \bar{i} \smallsetminus U_{j}}\p^{\alpha}_x\left[ \left\{ X_{t}^{i,1} \in \bar{i} \right\}\right]$. Fact \ref{fact2} is basic, to get it we use a similar method that \cite{Pierre1} use to get the Proposition 2.8 page 888. We now give a proof of Fact \ref{fact1}, 
\begin{eqnarray}
\pa\left[\bigcap_{i=1}^{n(f)}\bigcap_{x\in 
\bar{i} \smallsetminus U_{i}}\bigcap_{k=1}^{\eta(x)} \left\{ X_{t}^{x,k} \in \bar{i} \right\}\right]&=&\p^{\alpha}\left[\bigcap_{i=1}^{n(f)}\bigcap_{x\in 
\bar{i} \smallsetminus U_{i}}\bigcup_{k_x=0}^{+\infty}\left\{\eta(x)=k_x,\bigcap_{k=1}^{k_x} \left\{ X_{t}^{x,k} \in \bar{i} \right\}\right\}\right],  \\ 
 & = & \prod_{i=1}^{n(f)}\prod_{x\in 
\bar{i} \smallsetminus U_{i}}\p^{\alpha}_x\left[ \bigcup_{k_x=0}^{+\infty}\left\{\eta(x)=k_x,\bigcap_{k=1}^{k_x} \left\{ X_{t}^{x,k} \in \bar{i} \right\}\right\}\right], \nonumber \\
 & = & \prod_{i=1}^{n(f)}\prod_{x\in 
\bar{i} \smallsetminus U_{i}}\exp\left(\lambda(\p^{\alpha,1}_x\left[ \left\{ X_{t}^{x,1} \in \bar{i} \right\}\right]-1)\right) \nonumber .
 \end{eqnarray}
In the first equality we have introduced a partition over the values of the $(\eta(x),x\in \Z)$, the second equality comes from the independence of the particles starting from distinct points of the lattice,     the last one from the fact that $(X_.^{.,s}, s \in \N)$ are i.i.d. Finally
 \begin{eqnarray}
\pa\left[\bigcap_{i=1}^{n(f)}\bigcap_{x\in 
\bar{i} \smallsetminus U_{i}}\bigcap_{k=1}^{\eta(x)} \left\{ X_{t}^{x,k} \in \bar{i} \right\}\right] & \geq &\prod_{i=1}^{n(f)}\prod_{y\in 
\bar{i} \smallsetminus U_{i}}\exp\left(\lambda(a-1)\right) , \forall \alpha,
\end{eqnarray} 
so we get \ref{fact1}. To end the proof we use \ref{pr2} and \ref{pr3}. Notice that Fact 1 is true for all environments $\alpha$ whereas Fact 2 is only true for all good environments.  \\
To get \ref{3.1eq1} we will consider three different distances from the starting point of the particles to the boundary of $V_f$: $M_1$ or $M_{n(f)+1}$ (Cases 1-3 below). We will only discuss the case when the boundary is $M_{n(f)+1}$, the other case can be treated in the same way. \\
Case 1 (very long distance) notice that $\forall x > M_{n(f)+1}+t$, $X_{t}^{x,.}\notin V_f$, \\
Case 2 (long distance), let $ x \in J_t \equiv [M_{n(f)+1}+\log_2t(\log t)^2,M_{n(f)+1}+t]$. With the same computations we have obtained \ref{fact1} we can prove that
\begin{eqnarray}
\pa\left[\bigcap_{x \in J_t}  \bigcap_{k=1}^{\eta(x)}\{{X_{t}^{x,k} \notin V_{f}}\} \right] \equiv \pa\left[\bigcap_{x \in J_t}  \bigcap_{k=1}^{\eta(x)}\{{X_{t}^{x,k}> M_{n(f)+1}}\} \right] \geq    \exp\left(t\lambda(b-1)\right) , \forall \alpha
 \end{eqnarray}
 where  $b= \inf_{x \in J_t}\p^{\alpha}_x\left[ \left\{ X_{t} > M_{n(f)+1} \right\}\right]$, moreover it is a basic fact (see \cite{Pierre1} Proposition 2.8 page 888) that for all $\alpha \in G_t$, $b\geq 1- cte /(t (\log t)^{\gamma-2}).$\\  
 Case 3 (short distance), let $ x \in K_t \equiv [M_{n(f)+1}+(\log_2t)^2,M_{n(f)+1}+(\log_2t)(\log t)^2]$, we have  
 \begin{eqnarray}
\pa\left[\bigcap_{x \in K_t}  \bigcap_{k=1}^{\eta(x)}\{{X_{t}^{x,k} \notin V_{f}}\} \right] \geq    \exp\left((\log_2t)(\log t)^2\lambda(c-1)\right) , \forall \alpha,
 \end{eqnarray}
 where  $c= \inf_{x \in K_t}\p^{\alpha,1}_x\left[ \left\{ X_{t} > M_{n(f)+1} \right\}\right]$, moreover it is also a basic fact, see the reference above, that for all $\alpha \in G_t$, $c\geq 1- cte/(\log t)^{\gamma-2}$.\\  
From this three cases we deduce that the only particles, starting from a point inside $(M_{n(f)}+1,+\infty)$, that can reach $M_{n(f)}+1$, are the particles that belong to the interval $(M_{n(f)+1},M_{n(f)+1}+(\log_2t)^2)$, which is a subset of the set of indeterminate points $U_{M_n(f)+1}$, of size $(\log_2t)^2$, negligible comparing to $(\log  n)^2$.
\end{Pre}

 \begin{Pro} \label{3.2} For all $f\in C^1_{\kappa}$, $\gamma>12$, and all $t$ large enough, $Q[G_{t}] = 1-o(1)$ and 
 \begin{eqnarray}
\inf_{\alpha \in G_{t}}\left\{\pa\left[\bigcap_{i=1} ^{n(f)}\bigcap_{x\in 
\bar{i} \smallsetminus U_{i}}\bigcap_{k=1}^{\eta(x)} \left\{\left| \frac{X_{t}^{x,k}}{(\log t)^2}-m_{i}(t) \right| \leq \frac{cte}{(\log t)^{1/2}} \right\}\right]\right\} = 1-o(1).
 \end{eqnarray}
 \end{Pro}
 
\noindent The Proposition above generalizes Sinai's localization result for independent particles. It shows that each particle from a given valley will be located at time $t$ in a small neighborhood of the coordinate of the bottom of this same valley.  
\noindent \\

\begin{Pre}
The proof of this Proposition is based on the two following facts, let $I_{t,i}=(m_{i}(t) -cte (\log t)^{3/2},m_{i}(t) +cte (\log t)^{3/2})$, we have
\begin{eqnarray}
\pa\left[\bigcap_{i=1} ^{n(f)}\bigcap_{x\in 
\bar{i} \smallsetminus U_{i}}\bigcap_{k=1}^{\eta(x)} \left\{X_{t}^{x,k}\in I_{t,i} \right\}\right] \geq \exp\left(\lambda \left[d-1\right] \sum_{k=1}^{n(f)}|\bar{k}\smallsetminus U_k|\right) \nonumber
\end{eqnarray}
where $d=\inf_{1 \leq i \leq n(f)}\inf _{x \in \bar{i}\smallsetminus U_i} \p_i^{\alpha,1}\left[ X_t^{x,1} \in I_{t,i} \right]$, which is obtained with the same computation we did to get the previous Proposition, and 
\begin{eqnarray}
d \geq 1-\frac{cte}{(\log n)^{\gamma-12}}
\end{eqnarray}
which is obtained by using a similar method of \cite{Pierre1} (Theorem 2.11 page 889). 
\end{Pre}

\section{Proof of Theorem \ref{th1}}

In all this section we will always assume that $\alpha \in G_{t}$. 
The two Propositions \ref{3.1} and \ref{3.2} will be used frequently, as these results are true in $\pa$ probability we will mention it by using ''$\dot{=}$'' instead of the common ''$=$''. 
\noindent
First using the definition of $\eta(x,t)$ and \ref{3.1eq1}, we get 
\begin{eqnarray}
\sum_{y \in \Z} \eta(y,t) f(y/(\log t)^2) & \equiv &\sum_{x \in \Z} \sum_{i=1}^{\eta(x)} \sum_{y\in \Z}\un_{\{X_{t}^{x,i}=y\}} f(y/(\log t)^2) \\
& \peq& \sum_{l=1}^{n(f)} \sum_{x\in \bar{l}} \sum_{k=1}^{\eta(x)} \sum_{y\in \Z}\un_{\{X_{t}^{x,k}=y\}}  f(y/(\log t)^2) +o(\log^2 t).
\end{eqnarray}
Extracting the indeterminate particles and using \ref{3.1eq2}, we have
\begin{eqnarray}
\sum_{l=1}^{n(f)} \sum_{x\in \bar{l}} \sum_{k=1}^{\eta(x)} \sum_{y\in \Z}\un_{\{X_{t}^{x,k}=y\}}  f(y/(\log t)^2) & \equiv & \sum_{l=1}^{n(f)} \sum_{x\in \bar{l}-U_l} \sum_{k=1}^{\eta(x)}  \sum_{y \in \Z} \un_{\{X_{t}^{x,k}=y\}}  f(y/(\log t)^2)+\epsilon^{(1)}_n  \nonumber \\
& \peq &  \sum_{l=1}^{n(f)} \sum_{x\in \bar{l}-U_l} \sum_{k=1}^{\eta(x)}  \sum_{y \in \bar{l}} \un_{\{X_{t}^{x,k}=y\}}  f(y/(\log t)^2) +\epsilon^{(1)}_n,
\end{eqnarray}
where $\epsilon^{(1)}_n=\sum_{l=1}^{n(f)} \sum_{x\in U_l} \sum_{k=1}^{\eta(x)}  \sum_{y \in \Z} \un_{\{X_{t}^{x,k}=y\}}  f(y/(\log t)^2)$. Thanks to the fact that $f \in C_{\kappa}^1$, the law of large number for the sum of $\eta(x)$ and that by definition  $|U_.|=2(\log_2 t)^2$, it is easy to see that in $P_1$ probability $\epsilon^{(1)}_n=O(n(f) (\log_2 t)^2)$.
Proposition \ref{3.2} yields 
\begin{eqnarray}
\sum_{l=1}^{n(f)} \sum_{x\in \bar{l}-U_l} \sum_{k=1}^{\eta(x)}  \sum_{y \in \bar{l}} \un_{\{X_{t}^{x,k}=y\}}  f(y/(\log t)^2) & \peq & \sum_{l=1}^{n(f)} \sum_{x\in \bar{l}-U_{l}} \sum_{k=1}^{\eta(x)}   \sum_{y \in I(t,l)} \un_{\{X_{t}^{x,k}=y\}}  f(y/(\log t)^2)
\end{eqnarray}
where $I(t,l)=\{m_{l}(t) -cte (\log t)^{3/2},...,m_{l}(t) +cte (\log t)^{3/2}\}$. Now using the fact that $f\in C^1_{\kappa}$, we get that
\begin{eqnarray}
 & & \sum_{l=1}^{n(f)} \sum_{x\in \bar{l}-U_{l}} \sum_{k=1}^{\eta(x)}   \sum_{y \in I(t,l)} \un_{\{X_{t}^{x,k}=y\}}  f(y/(\log t)^2) \\
 &=& \sum_{l=1}^{n(f)} \left(f(m_{l}/(\log t)^2)+O\left(1/(\log t)^{1/2}\right)\right) \sum_{x\in \bar{l}-U_{l}} \sum_{k=1}^{\eta(x)}   \sum_{y \in I(t,f)} \un_{\{X_{t}^{x,k}=y\}}. \nonumber
\end{eqnarray}
Using once again Proposition \ref{3.2}, we get
\begin{eqnarray}
\sum_{l=1}^{n(f)} \sum_{x\in \bar{l}-U_{l}} \sum_{k=1}^{\eta(x)}   \sum_{y \in I(t,l)} \un_{\{X_{t}^{x,k}=y\}}  f(y/(\log t)^2)  \peq \sum_{l=1}^{n(f)}  \left(f(m_{l}/(\log t)^2)+O\left(1/(\log t)^{1/2}\right)\right)  \sum_{x\in \bar{l}-U_{l}}  \eta(x), \nonumber
\end{eqnarray}
collecting what we did above and using once again the law of large number, we get
\begin{eqnarray}
\sum_{y \in \Z} \eta(y,t) f(y/(\log t)^2) \peq \lambda  \sum_{l=1}^{n(f)}  \left(f(m_{l}/(\log t)^2)+O\left(1/(\log t)^{1/2}\right)\right)\left| \bar{l}-U_{l}\right| +o(n(f)\log^2 t).  \nonumber   
\end{eqnarray}
To end the proof we use \ref{pr1}, the fact that by definition $|U_.|=2(\log_2 t)^2$, and we divide what we get by $(\log t)^2$.

\section{Proof of Propositions \ref{ch1} and \ref{proenv} }

\textit{Proof of Proposition \ref{ch1}}
Using the vocabulary of  \cite{NevPit} and [1989b], the points $M_{i}$, and $m_{i}$ are $\Gamma_{t}$-extrema for a random walk. Moreover thanks to the theorem of Donsker (see for example \cite{Durrett1}, page 406) $(S_{tl}/\sqrt{l},l)$ converge in distribution to a two-sided Brownian motion, therefore we get that the points $L_i$ are 1-minima and $S_i$ 1-maxima for the two sided Brownian motion.  The first part of our Proposition \ref{ch1} follows from the Proposition page 241 of \cite{NevPit}, indeed the 1-extrema built a stationary renewal process, such that the difference of two consecutive extrema form a i.i.d. sequence with Laplace transform equal to $1/\cosh(\sqrt( {2} \lambda))$. Moreover $n(f)$ is the number of renewal within the support of $f$. It is easy to check that $\p[n(f)\leq k]=\p[R_k \leq |supp(f)|]$, where $R_k=\sum_{l=1}^k x_l$ with $x_l=L_{l+1}-L_l$ so we get that if $supp(f)$ is large the distribution of $n(f)$ is given by a normal law with mean $\sigma^2 supp(f)/2$ and variance $\sigma^4 supp(f)*3/4$.   \\
\textit{Sketch of the proof of Proposition \ref{ch1}} The first 3 properties can easily be deduced from the proof above. About Sinai's basic properties the proof can be found in the Appendix A of \cite{Pierre1}.

\noindent \\
\textbf{Acknowledgment} I would like to thank Errico Presutti for introducing me the subject and for several very helpful discussions. 

{\small \bibliography{thbiblio}}

\begin{thebibliography}{8}
\providecommand{\natexlab}[1]{#1}
\providecommand{\url}[1]{\texttt{#1}}
\expandafter\ifx\csname urlstyle\endcsname\relax
  \providecommand{\doi}[1]{doi: #1}\else
  \providecommand{\doi}{doi: \begingroup \urlstyle{rm}\Url}\fi

\bibitem[Ligget(1985)]{Ligget}
T.M. Ligget.
\newblock \emph{Interacting Particle systems}.
\newblock Springer Verlag, New York, 1985.

\bibitem[DeMasi and Presutti(1991)]{DeMPre}
A.~DeMasi and E.~Presutti.
\newblock \emph{Mathematical Methods for Hydrodynamic Limits}.
\newblock Lectures notes in mathematics, Springer Verlag, 1991.

\bibitem[Kipnis and Landim(1998)]{KipLan}
C.~Kipnis and C.~Landim.
\newblock \emph{Scaling limits of interracting particle systems}.
\newblock Springer, 1998.

\bibitem[Solomon(1975)]{Solomon}
F.~Solomon.
\newblock Random walks in random environment.
\newblock \emph{Ann. Probab.}, \textbf{3}\penalty0 (1):\penalty0 \ 1--31, 1975.

\bibitem[Sinai(1982)]{Sinai}
Ya.~G. Sinai.
\newblock The limit behaviour of a one-dimensional random walk in a random
  medium.
\newblock \emph{Theory Probab. Appl.}, \textbf{27}\penalty0 (2):\penalty0 \
  256--268, 1982.

\bibitem[Neveu and Pitman(1989)]{NevPit}
J.~Neveu and J.~Pitman.
\newblock \emph{Renewal property of the extrema and tree property of the
  excursion of a one-dimensional Brownian motion. S\'eminaire de
  Probabilit\'ees XXIII}, volume 1372.
\newblock Springer, 1989.

\bibitem[Andreoletti(2005)]{Pierre1}
P.~Andreoletti.
\newblock Alternative proof for the localisation of {S}inai's walk.
\newblock \emph{Journal of Statistical Physics}, 118:\penalty0 883--933, 2005.

\bibitem[Durrett(1996)]{Durrett1}
R.~Durrett.
\newblock \emph{Probability: Theory an Examples, 2nd ed.}
\newblock Duxbury Press, 1996.

\end{thebibliography}


\end{document}